\documentclass{amsart}
\usepackage{amssymb}

\newtheorem{thm}{Theorem}[section]
\newtheorem{cor}[thm]{Corollary}

\theoremstyle{definition}

\theoremstyle{remark}

\numberwithin{equation}{section}

\newcommand{\GS}{Gr\"obner-Shirshov}

\begin{document}

\title{Gr\"obner-Shirshov basis for the braid group in the Artin-Garside generators  }

\author[Bokut]{L. A. Bokut$^*$}
\address{Sobolev Institute of Mathematics, Novosibirsk 630090, Russia}

\email{bokut@math.nsc.ru}
\thanks{$^*$Supported in part by the Russian Fund of
    Basic Research, 05-01-00230, and the Integration Grant of the SB of the RAS, 1.9}

 \begin{abstract}
In this paper, we give a Gr\"obner-Shirshov basis of the braid group
$B_{n+1}$ in the Artin--Garside generators. As results, we obtain a
new algorithm for getting the Garside normal form, and a new proof
that the braid semigroup $B^+{n+1}$ is the subsemigroup in
$B_{n+1}$.
 \end{abstract}

\maketitle

\section{Introduction and the main theorem}

Markov \cite{Ma45} and Artin \cite{Ar47} independently found  a
normal form for words in the braid group

$$B_{n+1}=gp\langle a_1,\dots a_n|a_{i+1}a_ia_{i+1}=a_ia_{i+1}a_i, 1\leq i\leq n,
a_ka_s=a_sa_k, k-s>1
 \rangle.$$

It was proved in \cite{BCS} that the Markov-Artin normal form leads
to a \GS\ basis of $B_{n+1}$ in the Artin-Burau generators with the
so-called inverse tower order of words in the generators. Recall
that the Artin-Burau generators are the  elements
$$a_i, \\ A_{ij}=a_{j-1}\dots a_{i+1}a_i^2a_{i+1}^{-1}\dots a_{j-1}^{-1},$$

where $1\leq i <j \leq n+1$.

In the paper \cite{BFKS}, a \GS\ basis of the semigroup of positive braids $B_{n+1}^+$
in the Artin generators $a_i$ was found.

In this paper, we find a \GS\ basis of the braid group $B_{n+1}$ in the Artin-Garside generators
$a_i, 1\leq i \leq n, \Delta, \Delta^{-1}$ (\cite{Ga69}). Here we have

$$
\Delta=\Lambda_1\Lambda_2\dots \Lambda_n, \text{ \ with \ }
\Lambda_i=a_i\dots a_1.
$$
 Let us order these generators

$$
\Delta^{-1} < \Delta < a_1<\dots a_n.
$$

We order words in this alphabet in the deg-lex way  comparing two words first by theirs degrees
(lengths) and then lexicographically when the degrees are equal.

By $V(j,i), W(j,i), \dots$, where $j\leq i$, we  understand positive
words in the letters $a_j, a_{j+1}, \dots , a_i$. Also $V(i+1,i)=1,
W(i+1,i)=1,   \dots$.

Given $V=V(1,i)$, let $V^{(k)}, 1\leq k\leq n-i$ be the result of
shifting in $V$ all indices of all letters by $k$, $a_1\mapsto
a_{k+1}, \dots a_i\mapsto a_{k+i}$, and we also use the notation
$V^{(1)} =V'$. We write also $a_{ij}=a_ia_{i-1}\dots a_j, j\leq i-1,
a_{ii}=a_i, a_{ii+1}=1$.

\begin{thm}
A \GS\ basis of $B_{n+1}$ in the Artin-Garside generators consists of the following
relations :
\begin{align}
&a_{i+1}a_iV(1,i-1)W(j,i)a_{i+1j}=a_ia_{i+1}a_iV(1,i-1)a_{ij}W(j,i)',\label{E1}\\
 &a_sa_k=a_ka_s, s-k\geq 2,\label{E2}\\
&a_1V_1a_2a_1V_2\dots V_{n-1}a_n\dots a_1=\Delta V_1^{(n-1)}V_2^{(n-2)}\dots V_{n-1}',
\label{E3}\\
&a_l\Delta=\Delta a_{n-l+1},1\leq l\leq n,\label{E4}\\
& a_l\Delta^{-1}= \Delta^{-1} a_{n-l+1}, 1\leq l\leq n,\label{E4'}\\
&\Delta\Delta^{-1}=1, \Delta^{-1}\Delta=1,\label{E5}
\end{align}
where $1\leq i\leq n-1,$ $1\leq j\leq i+1, $
 $W$ begins with $a_i$ if it is not empty, and $V_i=V_i(1,i)$.
\end{thm}

Recall that a subset $S$ of the free algebra $k\langle X \rangle$
over a field $k$ on $X$ is called a \GS\ set (basis) if every
composition of elements of $S$ is trivial. This definition goes back
to  Shirshov's 1962 paper \cite{Sh62}. We recall the definition of
triviality of a composition in the below.

Let us define
$$
\Lambda_i^{(-)}=a_i\dots a_2, i\geq 2,\ \Lambda_1^{(-)}=1,
E_i=\Lambda_1\dots
\Lambda_{n-i}\Lambda_{n-i+1}^{(-)}\Lambda_{n-i+2}\dots \Lambda_n,
1\leq i\leq n.
$$

Then $E_ia_i=\Delta,$ and so \ $a_i^{-1}= \Delta^{-1}E_i.$ It
follows that we do not need the letters $a_i^{-1}$ and the relations
$a_ia_i^{-1}=1$, $a_i^{-1}a_i=1$ in the above presentation of the
group $B_{n+1}$.

We will write
$$
\Lambda_i^{(--)}=a_i\dots a_3, i\geq 3,\ \Lambda_2^{(--)}=1.
$$

\section{Proof of the Theorem }

Formulas (\ref{E1})-(\ref{E4}) are valid in $B_{n+1}^+$ for
\begin{align*}
&W(j,i)a_{i+1j}=a_{i+1j}W', V_{i-1}\Lambda_i=\Lambda_iV_{i-1}', a_i\Lambda_{i-1}\Lambda_i=
\Lambda_{i-1}\Lambda_ia_1,\\
& a_1\Lambda_{i+1}\dots \Lambda_n= \Lambda_{i+1}\dots
\Lambda_na_{n-i+1}.
\end{align*}

Formula (\ref{E5}) follows from (\ref{E4}).

Here and after notations are the same as in Theorem 1.1.

We need to prove that all compositions of relations
(\ref{E1})--(\ref{E5}) are trivial. The triviality of compositions
of (\ref{E1}), (\ref{E2}) was proved in \cite{BFKS}.

By "a word" we will mean a positive word in $a_i, \Delta$;
 $u=v$ is either the equality in $B_{n+1}^+$ or the letter-by-letter equality
(the meaning would be clear from the context).

 We use the following notation for words $u, v$:

$$
u\equiv v,
$$
if $u$ can be transformed to $v$ by the eliminations of leading
words of relations (\ref{E1})--(\ref{E4}), i.e., by the eliminations of
left parts of these relations. Actually, we will use an expansion of
this notation meaning that $u\equiv v$ if
$$
u\mapsto u_1\mapsto u_2\mapsto \dots \mapsto u_k= v,
$$
where $u_i< u$ for all $i$ and each transformation is an application
of (\ref{E1})--(\ref{E4}) (so, in general, only the first transformation
$u\mapsto u_1$ is the elimination
of the leading word of (\ref{E1})--(\ref{E4})).

An other expansion of that formula is
$$u\equiv v(mod\ w)$$
 meaning that $u$ can be transformed to $v$ as
before and all $u_i<w, u\leq w$.

By abuse of notations, we
take that in a word equivalence chain starting with a word $u$,
$$
u\equiv v\equiv w\equiv t \dots
$$
each equivalence $v\equiv w, \ w\equiv t, \dots $ is $mod\ u$.

This agrees with the definition of
triviality of a composition (see \cite{Bo72}, \cite{Bo76}). Namely,
a composition $(f,g)_w$ is called trivial $mod (S,w) $,  if
$$
(f,g)_w=\sum \alpha_ia_is_ib_i, a_i\overline{s_i}b_i<w,
$$
where $s_i\in S, a_i,b_i\in X^*, \alpha_i\in k $. Here $k\langle X
\rangle$ is a free associative algebra over a field $k$ on a set
$X$, $S\subset k\langle X \rangle$, $X^*$ is the set of all words in
$X$, $\overline{s}$ is the leading monomial of a polynomial $s$.
Recall that
$$
(f,g)_w=fb-ag, w=\overline{f}b=a\overline{g}, deg(f)+deg(g)>deg (w),
$$
or
$$
(f,g)_w=f-agb, w=\overline{f}=a\overline{g}b.
$$

Here $w$ is called the ambiguity of the composition $(f,g)_w$, $a,
b\in X^*$.

Let $S$ be the set of polynomial corresponding to semigroup
relations $u_i=v_i, u_i>v_i $, $(f,g)_w=u-v$ is a composition, $
u,v\in X^*, u,v<w$. The triviality  of $(f,g)_w$ $mod(S,w)$ means
that in the previous sense
$$
u\equiv t (mod\ w), v\equiv t (mod\ w)
$$
for some word $t$.

Take $V=V(j,i), 2\leq j\leq  i$ . By $V^{(-k)}$, where $1\leq k\leq
j-1$, we mean the
 result of shifting in $V$
the indices of all letters by $-k$, $a_j\mapsto a_{j-k}, \dots,
a_i\mapsto a_{i-k}$.

Take $\Delta_i=\Lambda_1\dots \Lambda_i, V=V(1,i), V^{\Delta_i^{\pm1}}=
\Delta_i^{\mp1} V \Delta_i^{\pm1}$. Then $ V^{\Delta_i^{\pm1}}$ is equal in
$B_{n+1}$ to
the  word that is the result of substitutions $a_j\mapsto a_{i-j+1}, 1\leq j\leq i$
in $V$. By abuse of notation, we will identify $ V^{\Delta_i^{\pm1}}$ with this word.

We need the formulas:

\begin{align}
&\Lambda_{i}W(2,i)\equiv W^{(-1)}\Lambda_{i},\
\Lambda_{i}^{(-)}W(3,i)
\equiv W^{(-1)}\Lambda_{i}^{(-)}\label{E6}\\
&a_i\Lambda_{i-1}\Lambda_i^{(-)}\equiv \Lambda_{i-1}\Lambda_i,\
a_i\Lambda_{i-1}\Lambda_i\equiv \Lambda_{i-1}\Lambda_ia_1,\label{E7}\\
&a_i\Lambda_{i-1}\Lambda_i^{(--)}\equiv \Lambda_{i-1}^{(-)}\Lambda_i,\label{E8}\\
&a_i\Lambda_{i-1}V_{i-1}\Lambda_i\equiv \Lambda_{i-1}\Lambda_ia_1V_{i-1}',\label{E9}\\
&a_i\Lambda_1V_1\dots \Lambda_{i-1}V_{i-1}\Lambda_i\equiv \Delta_ia_1V_1^{(i-1)}\dots
V_{i-1}',\label{E10}\\
&a_iV(1,i-1)\Lambda_1V_1\dots \Lambda_{i-1}V_{i-1}\Lambda_i
\equiv\Delta_ia_1V^{\Delta_i}V_1^{(i-1)}\dots V_{i-1}',\label{E11}\\
&a_iV(1,i-2)\Lambda_{i-1}W(2,i-1)\Lambda_i^{(-)}\equiv
\Lambda_{i-1}\Lambda_iV^{(2)}W',\label{E12}\\
&W'(1,i-1)\Lambda_2^{(-)}\dots \Lambda_i^{(-)}=
 \Lambda_2^{(-)}\dots
\Lambda_i^{(-)} W^{\Delta_i}, W=W(1,i-1), \label{E13}\\
&W(1,i-2)^{\Delta_{i-1}}\Lambda_i=\Lambda_iW^{\Delta_i},W=W(1,i-2),
\label{E13'}\\
 &a_n\cdot a_{n-1}a_n\cdot \cdot \cdot a_1\dots a_n\equiv
\Delta, \label{E14}
\end{align}
where  $i\geq 2, V_j=V_j(1,j), 1\leq j\leq i$.

Formula (\ref{E6}) is clear.

Formula (\ref{E7}) can be proved by induction on $i\geq 2$.
For $i=2$, it is clear. Let $i>2$.
Then
$$
a_i\Lambda_{i-1}\Lambda_i^{(-)}=a_ia_{i-1}\Lambda_{i-2}a_i\Lambda_{i-1}^{(-)}\equiv
a_{i-1}a_ia_{i-1}\Lambda_{i-2}\Lambda_{i-1}^{(-)}\equiv\\
a_{i-1}a_i\Lambda_{i-2}\Lambda_{i-1}\equiv \Lambda_{i-1}\Lambda_i.
$$

Formula (\ref{E8}) can also be proved by induction on $i\geq 2$. For $i=2$, it is clear. Let $i>2$.
Then
$$
a_i\Lambda_{i-1}\Lambda_i^{(--)}=a_ia_{i-1}\Lambda_{i-2}a_i\Lambda_{i-1}^{(--)}\equiv
a_{i-1}a_ia_{i-1}\Lambda_{i-2}\Lambda_{i-1}^{(--)}\equiv\\
a_{i-1}a_i\Lambda_{i-2}^{(-)}\Lambda_{i-1}\equiv \Lambda_{i-1}^{(-)}\Lambda_i.
$$
Formula (\ref{E9}) follows from (\ref{E1}) and (\ref{E6}).

Formula (\ref{E10}) is clear for $i=2$; for $i>2$, it follows from
(\ref{E1}) and (\ref{E7}):
\begin{align*}
&a_i\Lambda_1V_1\dots \Lambda_{i-1}V_{i-1}\Lambda_i\equiv
\Lambda_1V_1\dots \Lambda_{i-2}V_{i-2}a_ia_{i-1}\Lambda_{i-2}V_{i-1}\Lambda_i\equiv\\
&\Lambda_1V_1\dots
\Lambda_{i-2}V_{i-2}a_{i-1}a_ia_{i-1}\Lambda_{i-2}\Lambda_{i-1}V_{i-1}'\equiv
\Lambda_1V_1\dots
\Lambda_{i-2}V_{i-2}a_{i-1}a_i\Lambda_{i-2}\Lambda_{i-1}a_1V_{i-1}'\equiv\\
&\Lambda_1V_1\dots
\Lambda_{i-2}V_{i-2}\Lambda_{i-1}\Lambda_ia_1V_{i-1}'\equiv
\Delta_ia_1V_1^{(i-1)}\dots V_{i-1}'.
\end{align*}

Remark, that the last word is less then the first word in the above
chain of equivalence formulas, though  it can be greater then the
word just before the last.

Formula (\ref{E11}) can be proved by induction on $i\geq 2$. It is
clear for $i=2$. Take $i>2$. Then formula (\ref{E11}) follows from
(\ref{E10}) by induction on the number $k$ of letters $a_{i-1}$ in
$V(1,i-1)$. If $k=0$, the result is clear. Take $k\geq 1$ and
$$V(1,i-1)=W(1,i-1)a_{i-1}T(1,i-2).$$

Then
\begin{align*}
&a_{i-1}T(1,i-2)\Lambda_1V_1\dots \Lambda_{i-2}V_{i-2}\Lambda_{i-1}
\equiv\Delta_{i-1}a_1T^{\Delta_{i-1}}V_1^{(i-2)}\dots V_{i-2}'.
\end{align*}

It follows that
\begin{align*}
&a_iW(1,i-1)a_{i-1}T(1,i-2)\Lambda_1V_1\dots \Lambda_{i-1}V_{i-1}\Lambda_i\equiv \\
&a_iW(1,i-1)\Delta_{i-1}a_1T^{\Delta_{i-1}}V_1^{(i-2)}\dots V_{i-2}'V_{i-1}\Lambda_i\equiv\\
&a_iW(1,i-1)\Delta_ia_2T^{\Delta_i}V_1^{(i-1)}\dots V_{i-1}'\equiv
\Delta_ia_1W^{\Delta_i}a_2 T^{\Delta_i}V_1^{(i-1)}\dots V_{i-1}'.
\end{align*}

Remark again that here the second and third words above are less then
the first one for $i-1\geq 2$.

Formula (\ref{E14}) can be proved by induction on $n$. Indeed:

\begin{align*}
&a_n\cdot a_{n-1}a_n\cdot \cdot \cdot a_2\dots a_n\cdot  a_1\dots a_n\equiv
a_2\cdot a_3a_2 \cdots a_n\dots a_2 \cdot a_1\dots a_n \equiv \\
&a_2\cdot a_3a_2 \cdots a_{n-1}\dots a_2 \cdot a_1\dots a_{n-1}\cdot a_n a_{n-1}\dots a_1 \equiv
\dots \equiv \\
&a_2\cdot a_3a_2 \cdot a_1a_2a_3 \cdot a_4\dots a_1 \cdots a_n\dots
a_1 \equiv a_2\cdot a_1a_2\cdot a_3a_2a_1\cdots a_n\dots a_1\equiv
\Delta.
\end{align*}

Now let us check the composition of (\ref{E1}), $i+1=n$, and
(\ref{E3}). Without loss of generality we may assume that $j=1$.

The ambiguity is
$$
w=a_na_{n-1}V(1,n-2)W(1,n-1)\Lambda_1V_1\dots\Lambda_{n-1}V_{n-1}a_{n1}.
$$

Applying (\ref{E3}) to $w$, we obtain
$$
w_1=a_na_{n-1}V(1,n-2)W(1,n-1)\Delta V_1^{(n-1)}\dots V_{n-1}'.
$$

Applying (\ref{E1}) to $w$, we obtain
$$
w_2= a_{n-1}a_na_{n-1}V(1,n-2)\Lambda_{n-1}W'\Lambda_1'V_1'\dots
\Lambda_{n-1}'V_{n-1}'.
$$

Applying (\ref{E1}) and (\ref{E6})-(\ref{E14})  to $w_1, w_2$, we
obtain
\begin{align*}
&w_1\equiv \Delta a_1a_2V^{\Delta}W^{\Delta}V_1^{(n-1)}\dots V_{n-1}',\\
&w_2= a_{n-1}a_na_{n-1}V(1,n-2)\Lambda_{n-1}W'\Lambda_2^{(-)}V_1'\dots \Lambda_{n-1}^{(-)}
V_{n-2}'\Lambda_n^{(-)}V_{n-1}'\equiv \\
&a_{n-1}a_na_{n-1}V(1,n-2)\Lambda_{n-1}W'\Lambda_2^{(-)}\dots \Lambda_{n}^{(-)}V_1^{(n-1)}
\dots V_{n-1}' (mod\ w)\equiv \\
&a_{n-1}a_na_{n-1}V(1,n-2)\Lambda_{n-1}\Lambda_2^{(-)}\dots \Lambda_{n-1}^{(-)}\Lambda_{n}^{(-)}
W^{\Delta}V_1^{(n-1)}\dots V_{n-1}' (mod\ w)\equiv \\
&a_{n-1}a_na_{n-1}V(1,n-2)\Lambda_1\Lambda_2\dots \Lambda_{n-2}\Lambda_{n-1}\Lambda_{n}^{(-)}
W^{\Delta}V_1^{(n-1)}\dots V_{n-1}' (mod\ w)\equiv \\
&a_{n-1}a_na_{n-1}\Lambda_1\Lambda_2\dots \Lambda_{n-2}\Lambda_{n-1}\Lambda_{n}^{(-)}V^{\Delta}
W^{\Delta}V_1^{(n-1)}\dots V_{n-1}' (mod\ w)\equiv \\
&a_{n-1}a_n\Lambda_1\dots \Lambda_{n-2}\Lambda_{n-1}\Lambda_{n}^{(--)}a_1a_2V^{\Delta}
W^{\Delta}V_1^{(n-1)}\dots V_{n-1}' (mod\ w)\equiv \\
&a_{n-1}\Lambda_1\dots \Lambda_{n-2}\Lambda_{n-1}^{(-)}\Lambda_{n}
a_1a_2V^{\Delta}W^{\Delta}V_1^{(n-1)}\dots V_{n-1}' (mod\ w)\equiv \\
&\Delta a_1a_2V^{\Delta}W^{\Delta}V_1^{(n-1)}\dots V_{n-1}'(mod\ w).
\end{align*}
The composition is checked.

There is the composition of (\ref{E1}), $i+1<n$ and (\ref{E3}).
Again, we may assume that $j=1$ .

The ambiguity is
$$
w=a_{i+1}a_iV(1,i-1)W(1,i)\Lambda_1V_1\dots\Lambda_{i+1}V_{i+1}\dots V_{n-1}\Lambda_n.
$$

Applying (\ref{E3}) to $w$, we obtain
\begin{align*}
&w_1=a_{i+1}a_iV(1,i-1)W(1,i)\Delta V_1^{(n-1)}\dots V_{n-1}'\equiv \\
&\Delta a_{n-i}a_{n-i+1}V^{\Delta}W^{\Delta}V_1^{(n-1)}\dots V_{n-1}'
\end{align*}

Applying (\ref{E1}) to $w$, we have
\begin{align*}
&w_2= a_ia_{i+1}a_iV(1,i-1)\Lambda_iW(1,i)'\Lambda_1'V_1'\dots
\Lambda_i'V_i'V_{i+1}\dots V_{n-1}\Lambda_n \equiv \\
&\Lambda_1\dots \Lambda_{i+1}a_1a_2V^{\Delta_{i+1}}W^{\Delta_{i+1}}V_1^{(i)} \dots V_i'V_{i+1}
\Lambda_{i+2}\dots V_{n-1}\Lambda_n (mod\ w)\equiv \\
&\Delta a_{n-i}a_{n-i+1}V^{\Delta}W^{\Delta}
V_1^{(n-1)} \dots V_i^{(n-i)}\dots V_{n-1}'. \\
\end{align*}

Here we use the calculation of $w_2$ from the previous composition,
substituting $n\mapsto i+1$.

The second composition is checked.

There is the composition of (\ref{E3}) and (\ref{E1}).

The ambiguity is
$$
w=\Lambda_1V_1\dots V_{n-1}a_{ni+2}a_{i+1}a_i\dots a_1V(1,i-1)W(j,i)a_{i+1j}.
$$

Applying (\ref{E3}) to $w$, we  obtain
\begin{align*}
&w_1=\Delta V_1^{(n-1)}\dots V_{n-1}'VWa_{i+1j}.
\end{align*}

Applying (\ref{E1}) to $w$, we have
\begin{align*}
&w_2= \Lambda_1V_1\dots V_{n-1}a_{ni+2}a_ia_{i+1}a_i\dots a_1V(1,i-1)a_{ij}W'\equiv \\
&\Lambda_1V_1\dots V_{n-1}a_ia_{ni+2}a_{i+1}a_i\dots a_1Va_{ij}W'\equiv
\Delta V_1^{(n-1)}\dots V_{n-1}'a_{i+1}Va_{ij}W'\equiv \\
&\Delta V_1^{(n-1)}\dots V_{n-1}'Va_{i+1j}W'\equiv \Delta V_1^{(n-1)}\dots V_{n-1}'VWa_{i+1j}.
\end{align*}

There is the composition of (\ref{E1}) and (\ref{E4}).

The ambiguity is
$$
w=a_{i+1}a_iV(1,i-1)W(j,i)a_{i+1j+1}a_j\Delta.
$$

Applying (\ref{E4}) to $w$, we obtain
\begin{align*}
&w_1=a_{i+1}a_iVWa_{i+1j+1}\Delta a_{n-j+1} \equiv
\Delta a_{n-i}a_{n-i+1}V^{\Delta}W^{\Delta}a_{n-i}a_{n-i+1}\dots a_{n-j+1} \equiv \\
&\Delta a_{n-i}a_{n-i+1}V^{\Delta}a_{n-i}a_{n-i+1}\dots a_{n-j+1}W^{\Delta (-1)}\equiv \\
&\Delta a_{n-i}a_{n-i+1}a_{n-i}V^{\Delta}a_{n-i+1}\dots a_{n-j+1}W^{\Delta (-1)}.
\end{align*}

Applying (\ref{E1}) to $w$, we  obtain
\begin{align*}
&w_2=a_ia_{i+1}a_iVa_{ij}W'\Delta \equiv \Delta
a_{n-i+1}a_{n-i}a_{n-i+1}
V^{\Delta}a_{n-i+1}\dots a_{n-j+1}W'^{\Delta} \equiv \\
&\Delta a_{n-i}a_{n-i+1}a_{n-i}V^{\Delta}a_{n-i+1}\dots a_{n-j+1}W'^{\Delta}.
\end{align*}

The composition is checked for $W^{\Delta (-1)}=W'^{\Delta}$.

The triviality of the composition of (\ref{E1}) with  (\ref{E4'}) is
proved similarly.

There is the composition of (\ref{E3}) and (\ref{E4}).

The ambiguity is
$$
w=\Lambda_1V_1\dots V_{n-1}a_{n}\dots a_2a_1\Delta.
$$

Applying (\ref{E3}) to $w$, we  have
\begin{align*}
&w_1=\Delta V_1^{(n-1)}\dots V_{n-1}'\Delta\equiv \Delta^2V_1^{(n-1)\Delta}\dots V_{n-1}^{'\Delta}.
\end{align*}

Applying (\ref{E4}) to $w$, we obtain
\begin{align*}
&w_2= \Lambda_1V_1\dots V_{n-1}a_{n}\dots a_2\Delta a_n \equiv \Delta
\Lambda_1^{\Delta}V_1^{\Delta}\dots V_{n-1}^{\Delta}a_{1}\dots a_{n-1}a_n = \\
&\Delta a_nV_1^{\Delta}a_na_{n-1}\dots V_{n-1}^{\Delta}a_{1}\dots a_n \equiv
\Delta a_n\cdot a_{n-1}a_n\cdots a_1\dots a_n V_1^{\Delta (-(n-1))}\dots V_{n-1}^{\Delta
(-1)}\\
&\equiv \Delta^2V_1^{(n-1)\Delta}\dots V_{n-1}^{'\Delta}
\end{align*}

for $V_i^{(j)\Delta}=V_i^{\Delta (-j)}, 2\leq i+j\leq n$.

The composition is checked.

The triviality of composition of (\ref{E3}) and (\ref{E4'}) is
proved similarly.

\section{Corollaries}

Let $S\subset k\langle X \rangle$.
A word $u$ is called $S$-irreducible
if $u\neq a\overline{s}b$, where $s\in S, a,b \in X^*$. Let $Irr(S)$ be the set
of all $S$-irreducible words.

Recall  Shirshov's Composition lemma (\cite{Sh62}, \cite{Bo72},
\cite{Bo76}) :

{\it Let $S$ be a \GS\ set in $k\langle X \rangle$. If $f\in ideal(S)$, then
$\overline{f}=a\overline{s}b, s\in S$. The converse is also true.}

The Main corollary to this lemma is the following statement:

{\it Subset $S \subset k\langle X \rangle $ is a \GS\ set iff the
set $Irr(S)$ is a linear basis for the algebra $k\langle X \rangle /
ideal(S)= k\langle X| S \rangle$ generated by $X$ with defining
relations $S$.}

Let $G=sgp\langle X|S \rangle$ be the semigroup generated by $X$
with defining relations $S$. Then $S$ is called a \GS\ basis of $G$
if $S$ is a \GS\ basis of the semigroup algebra $k(G)$, i.e., $S$ is
a \GS\ set in $ k\langle X \rangle$. It follows from the Main
corollary to the Composition lemma that in this case any word $u$ in
$X$ is equal in $G$ to a unique $S$-irreducible word $C(u)$, called
the normal (canonical) form of $u$.

Now let $S$ be the set of relations (\ref{E1})-(\ref{E5}), and let
$C(u)$ be a normal form of a word $u\in B_{n+1}$. Then $C(u)$ has a
form
$$
C(u)=\Delta^k A,
$$
where $k\in \mathbb{Z}$, and  $A$ a positive $S$-irreducible word in
$a_i$'s. Let us prove that $A\neq \Delta A_1$ for every positive
word $A_1$.

Let us switch to the semigroup $B_{n+1}^+$. Note that we have no
generator $\Delta$ in $B_{n+1}^+$. The set $S_1$ of
(\ref{E1})-(\ref{E2}) is a \GS\ basis of $B_{n+1}^+$. Then $A$  is
an $S_1$-irreducible word.

Suppose that
$$
A= \Delta A_1,
$$
where $A_1$ is a positive word. We may assume that $A_1$ is an $S_1$-irreducible word too.  Then
$$
A=a_1\cdot a_2a_1\cdots a_n\dots a_1 A_1
$$
in the semigroup $B_{n+1}^+$. Let
us prove that
$$
A= C(a_1\cdot a_2a_1\cdots a_n\dots a_1 A_1)
=\Lambda_1V_1\dots \Lambda_{n-1}V_{n-1}\Lambda_n V,
$$
where here $=$ is the graphical equality,   $V_i=V(1,i), 1\leq i\leq
n-1$, and $C(D)$ is a normal form of $D$ in $B_{n+1}^+$.  It would
contradict the $S$-irreducibility of $A$.

 More generally, let us prove that
$$
 C(\Lambda_1W_1\dots \Lambda_{n-1}W_{n-1}\Lambda_n W_n)
=\Lambda_1V_1\dots \Lambda_{n-1}V_{n-1}\Lambda_n V_n,
$$
where $W_i=W_i(1,i), V_i=V_i(1,i), 1\leq i\leq n$. If
$$
B= \Lambda_1W_1\dots \Lambda_{n-1}W_{n-1}\Lambda_n W_n
$$
is a canonical word, than we are at home. Suppose, $B$ contains the
left part of a relation (\ref{E1})-(\ref{E2}). If this word is a
subword of $W_i, 1\leq i\leq n $, then the situation is clear: we
can apply the relation to get a smaller  word of the same form, and
then we can use induction. Let the word be a subword of
$\Lambda_kW_k,1\leq k\leq n$, but not $W_k$. It may only be  the
left part of (\ref{E1}), $k\geq i+1$. Then we have
\begin{align*}
&\Lambda_kW_k= a_{ki+2}a_{i+1}a_ia_{i-1}\dots a_1V(1,i-1)W(j,i)a_{i+1j}T_k\equiv\\
& a_{ki+2}a_ia_{i+1}a_ia_{i-1}\dots a_1V(1,i-1)a_{ij}W'T_k \equiv
a_i\Lambda_k W_{k1},
\end{align*}
where $T_k=T_k(1,k), W_{k1}=W_{k1}(1,k)$. Substituting this quantity
of $\Lambda_kW_k$ in $B$, we obtain a positive word $D$ which is
smaller then $B$ and has the same form. By  induction, we are at
home.

As a result, we have the following

\begin{cor} The $S$-irreducible normal form of each word of $B_{n+1}$ coincides
with the Garside normal form of the word.
\end{cor}
\begin{proof}Recall (\cite{Ga69}) that Garside normal form $G(u)$ of $u\in B_{n+1}$
is
$$
G(u)=\Delta^k A,
$$
where $u=G(u)$ in $B_{n+1}$, $k\in \mathbb{Z}$, and $A$ is a
positive word in $a_i$, $A\neq \Delta A_1$ for every positive word
$A_1$, and $A$ is the minimal word  with these properties. We have
proved that the $S$-irreducible normal form $C(u)$ has this
properties.
\end{proof}

\begin{cor}( \cite{Ga69} ) The semigroup of positive braids $B_{n+1}^+$ can be embedded  into a
group.
\end{cor}
\begin{proof}  It follows immediately from Theorem 1.1 that $B_{n+1}^+\subset B_{n+1}$.
\end{proof}

\end{document}